\documentclass{arSIH}

\newtheorem{theorem}{Theorem}
\newtheorem{proposition}[theorem]{Proposition}

\theoremstyle{definition}
\newtheorem{definition}[theorem]{\font\=cmssi10\Definition\bf}
\newtheorem{example}[theorem]{\font\=cmssi10\Example\bf}
\newtheorem{remark}[theorem]{\font\=cmssi10\Remark\bf}

\def\Cal{\mathcal}
\def\Eps{\hbox{\font\=cmmi10 scaled\magstep1\\char'017}\kern0.15mm}% bigger varepsilon
\def\Iota{\kern.15mm\hbox{\font\=cmmi10 scaled\magstep1\\char'023}\kern0.2mm}% bigger iota
\def\Nu{\hbox{\font\=cmmi10 scaled\magstep1\\char'027}\kern0.25mm}% bigger nu
\def\uvarPi{\kern.15mm\underline{\kern-.15mm\varPi\kern-.85mm}\kern.85mm}% underlined varPi
\def\uOmega{\kern.3mm\underline{\kern-.3mm\Omega\kern-.3mm}\kern.3mm}% underlined Omega

\def\fRe{\hbox{\font\=cmr9\f\kern.1mm}\roman{Re}\kern.75mm}% the real      part of a function
\def\fIm{\hbox{\font\=cmr9\f\kern.1mm}\roman{Im}\kern.65mm}% the imaginary part of a function
\def\vecc#1{\kern-.5mm\vec{\kern.5mm#1}}% \vec moved a bit backwords so that in (\vecc x) the arrow would not get too close to ')'
\def\TVS{\roman{TVS}\kern0.37mm}%
\def\LCS{\roman{LCS}\kern0.37mm}%
\def\BaS{\roman{BaS}\kern0.37mm}%
\def\dimHa{{\rm dim_{_{\kern.2mm Ha}}}}% the Hamel dimension of a structured vector space
\def\Linb{\Cal L\lower.7mm\hbox{\kern.1mm\font\=cmmi6\b}}% \Linb(E,F) gives \Cal L_b(E,F) = the topological vector space of continuous linear maps E to F with the topology of uniform convergence on bounded sets
\def\LL^#1{L\kern0.15mm\raise.4mm\hbox{$^{#1}$}\kern0.15mm}% = L^p with p a bit lifted

\def\dualbeta{^{\kern0.4mm\prime}_{\kern-.2mm\raise.95mm\hbox{$_{_\beta}$}}} % E\dualbeta = E'_\beta = strong topological dual of E
\def\Nbh{\Cal N_{\font\=cmmi6\lower.15mm\hbox{\kern.1mm\bh\kern.15mm}}}% \Nhb(x,T) = set of T-neighborhoods of x
\def\Topma{\roman{{Top_{}}_{\hbox{\font\=cmr6\ma}}\kern.15mm}}
\def\co{\hbox{\font\=cmmi12\c}\kern.15mm\lower.15mm\hbox{$_{\rm o}$}}% \co(I,F) = the supremumnormable space of functions I to F with become small outside finite sets
\def\prodc{\prod{_{_{\kern-.3mm\bold c\kern.15mm}}}}% cartesian propduct of sets
\def\vsprod_#1_#2{\prod\kern-0.3mm{}_{_{\roman{#1}\sp{#2}\,}}} % for example, \vsprod_TVS_\Re produces \prod_{\roman{TVS} I\!\!R}
\def\vscoprod_#1_#2{\coprod\kern-0.3mm{}_{_{\roman{#1}\sp{#2}\,}}} % see the preceding
\def\expnota^#1]_#2{\,^{#1\,]{_{}}_{\roman{#2}}}} % F\expnota^\Omega_{tvs} = the topological product vector space of all functions Omega to F

\def\bold#1{{\bf#1}}
\def\roman#1{{\rm#1}}
\def\limu_#1{\lim\kern-5.5mm\lower1.5mm\hbox{$_{#1}\ $}}
\def\oseoy{\raise1.9mm\hbox{\kern.5mm\font\=cmr5\o}\kern-1.7mm y}% kahdessa kohtaa esiintyy
\def\O{{}^{}\Cal O}
\def\Pows{\Cal P\kern-.4mm_s\kern.3mm}
\def\lei{      {}_{ {}^{\,\downarrow\text{\hskip-2.1mm}       }  }  \cap       }
\def\lei{\hbox{\kern.45mm$_{^\downarrow}\kern-1.280mm\cap\kern.85mm$}}
\def\Ze{Z\!\!\!Z} %                                                kokonaisluvut
\def\inve{\lower.85mm\hbox{$^{^-}$}\kern-.5mm{}^\iota}

\def\fvalue{\hbox{\kern.2mm\font\=cmr10\\char'022\kern-.2mm}} % f\value x=f(x)
\def\ffvalue{\hbox{\kern.2mm\font\=cmr7\\char'022\kern-.2mm}} % f\value x=f(x)
\def\image{\hbox{\font\=cmr10\\char'022\kern-1mm\char'022}} % f\image A=f[A]
\def\iimage{\hbox{\font\=cmr7\\kern.3mm\char'022\kern-.7mm\char'022\kern-.3mm}} % f\image A=f[A]
\def\images{\hbox{\font\=cmr10\\char'022\kern-1mm\char'022\kern-1mm\char'022}} % f\images\Cal A={f[A]:A in Cal A}
\def\timesn{\kern-.2mm\times\kern-.2mm} % same as times but with smaller spaces before and after
\def\ttimes{\hbox{\kern-.2mm${}\times\kern-2.5mm\lower.8mm\hbox{\font\=cmr5\t}\kern1.8mm$}} % product topology = \rist
\def\ktimes{\hbox{\kern-.2mm${}\times\kern-2.5mm\lower1mm\hbox{\font\=cmr5\k}\kern1.5mm$}} % compactly generated product topology
\def\vstimes{\kern.95mm\raise.45mm\hbox{\font\=cmbsy6\\char'002}\kern-2.3mm\lower.9mm\hbox{\font\=cmr5\vs}\kern1.05mm} % X\vstimes Y = the vector space product of X and Y
\def\Examplee{{\font\=cmssi10\E\kern.15mmx\kern.15mma\kern.15mmm\kern.14mmp\kern.17mml\kern.15mme}\kern.3mm. }
\def\Examples{{\font\=cmssi10\E\kern.15mmx\kern.15mma\kern.15mmm\kern.14mmp\kern.17mml\kern.15mme\kern.15mms}\kern.3mm. }
\def\Remarkk{{\font\=cmssi10\R\kern.15mme\kern.15mmm\kern.15mma\kern.15mmr\kern.15mmk\kern.15mm. }}
\def\Remarkss{{\font\=cmssi10\R\kern.15mme\kern.15mmm\kern.15mma\kern.15mmr\kern.15mmk\kern.15mms}\kern.3mm. }
\def\N{{I\!\!N}} % former \Ne, the numbers 1,2,...
\def\No{{I\!\!N\kern-.54mm\lower.15mm\hbox{$_{\rm o}$}}} % former \Neo, the natural numbers 0,1,2,...
\def\Nopot#1{I\!\!N\kern-.54mm\lower.15mm\hbox{$_{\rm o}$}\kern-.7mm{}^{#1}} % N_o^{#1}
\def\potNo{^{\kern.37mm I\!\!{N_{}}_{\kern-.22mm{\rm o}}}} % ^{IN_o}
\def\minus{\kern.2mm\lower1.05mm\hbox{$^-$}}
\def\pplus{\raise.22mm\hbox{\font\=cmr5\\char'053}}% 5 point +
\def\mminus{\raise.18mm\hbox{\font\=cmsy5\\char'000}}% 5 point -
\def\plusinftyy{\raise.18mm\hbox{\font\=cmr5\\char'053}\infty}% +infty for sub- and superscripts with smaller +
\def\minusinftyy{\raise.18mm\hbox{\font\=cmsy5\\char'000}\infty}% -infty for sub- and superscripts with smaller -
\def\plusinfty{\lower1.05mm\hbox{$^+$}\infty}
\def\minusinfty{\lower1.05mm\hbox{$^-$}\infty}
\def\Qe{\hbox{$Q\kern-2.6mm\raise.2mm\hbox{\font\=cmssqi8\I}\kern1.7mm$}}% the set of rational numbers
\def\Re{{I\!\!R}}
\def\Rep{{{I\!\!R^{\phantom{l}}}^{{}_{{}^{\!}\!+}}}}

\def\Ce{{\hbox{$C\kern-2.5mm\raise.2mm\hbox{\font\=cmssqi8\I}\kern1.48mm$}}}
\def\imag{\kern.15mm\lower.6mm\hbox{$^{^*}$}\kern-1.8mm\imath\kern.1mm} % the imaginary unit
\def\ebiF{\kern.1mm\hbox{\font\=cmmib8\F}\kern.5mm} % 8 point bold italic F
\def\ebiT{\kern.1mm\hbox{\font\=cmmib8\T}\kern.6mm} % 8 point bold italic T
\def\ebiU{\kern.1mm\hbox{\font\=cmmib8\U}\kern.5mm} % 8 point bold italic U
\def\fssi#1{\hbox{\font\=cmssi10\#1}\kern0.15mm} % cmssi in math
\def\smb#1{\hbox{\font\†=cmmi8\†#1\kern.3mm}} % eightpoint (capital) math symbols
\def\eCal#1{\kern.1mm\hbox{\font\†=cmbsy8\†#1\kern.4mm}} % 8point bold calligraphic #1
\def\ecal#1{\kern.1mm\hbox{\font\†=cmsy8\†#1\kern.3mm}} % 8point calligraphic #1
\def\ncal#1{\kern.1mm\hbox{\font\†=cmsy9\†#1\kern.3mm}} % 9point calligraphic #1
\def\vcal#1{\kern-.1mm\vec{\kern.2mm\hbox{\font\†=cmsy7\†#1}\kern.3mm}} % e.g., vector bundle \vcal E

\def\idv{\hbox{\font\=cmr10\id}\kern.25mm\lower.8mm\hbox{\font\=cmr7\v}\kern.3mm} % id_v E=id(v_s E)
\def\seq#1{\langle#1\rangle}
\def\Seq#1{\big\langle#1\big\rangle}

\def\SemiNor{\Cal S_{_N}\kern0.15mm}% \SemiNor E = set of continuous seminorms in E
\def\vecs{\upsilon\kern-0.3mm\lower.15mm\hbox{$_s$}\kern0.2mm} % underlying set of a structured vector space
\def\vecss{\hbox{\font\=cmitt10\v}\kern-0.1mm\lower.15mm\hbox{$_s$}\kern0.2mm} % underlying set of a vector space / vector (= linear) struc
\def\bnull#1{\hbox{\font\=cmssbx10\0}{}_{\font\=cmmi6\lower.15mm\hbox{\kern-.1mm\#1\kern.15mm}}} % \bnull X = zero in vector space X
\def\bzero#1{\hbox{\font\=cmbx10\0}{}_{\font\=cmmi6\lower.15mm\hbox{\kern-.1mm\#1\kern.15mm}}} % \bzero X = zero in vector space X
\def\dom{{{}^{}{\rm dom}\,{}_{{}^{}}}}
\def\domm{{}^{}{\rm dom}^{\kern.3mm\hbox{\font\=cmr6\2}}\,}

\def\rng{{}^{}{\rm rng}\,{}_{{}^{}}}
\def\CPi#1{C\kern-.2mm\lower.05mm\hbox{$_{_\Pi}$}\kern-1.52mm{}^{#1}}
\def\CinftyPi{C\kern.4mm\raise.3mm\hbox{$^\infty$}\kern-3.35mm_{_\Pi}\kern1.45mm}
\def\CinftyS{\Cinfty\kern-3.9mm_{_{\Cal S}}\kern1.45mm}
\def\Cinfty{C\kern.4mm\raise.3mm\hbox{$^\infty$}\kern.15mm}
\def\Cinftyzero{\hbox{$C\kern.4mm\raise.3mm\hbox{$^\infty$}\kern.15mm\kern-3.5mm_{\font\=cmr6\lower.15mm\hbox{\kern.1mm\0}}\kern1.9mm$}}

\def\RHB#1#2{\raise#1mm\hbox{$#2$}} % raised (by #1 mm) horizontal box of #2
\def\LHB#1#2{\lower#1mm\hbox{$#2$}} % lowered (by #1 mm) horizontal box of #2

\def\fiveroman#1{\hbox{\font\=cmr5\#1\kern.1mm}}

\def\sixroman#1{\hbox{\font\=cmr6\#1\kern.1mm}}
\def\eightmath#1{\hbox{\font\=cmmi8\{#1}\kern.1mm}}% 8 point math italic
\def\eightroman#1{\hbox{\font\=cmr8\{#1}\kern.1mm}}% 8 point roman
\def\subtext#1{\raise.2mm\hbox{$_{_{\kern0.15mm\roman{#1}}}$}}
\def\subtexT#1{\raise.2mm\hbox{$_{_{\kern0.15mm\hbox{\font\=cmr5\#1}}}$}}
\def\sNor#1{\kern.25mm\lower.38mm\hbox{$_{#1}$}}
\def\sNorr#1{\kern-.2mm\lower.38mm\hbox{$_{#1}$}}
\def\sNoreset_#1{\kern.13mm\lower.83mm\hbox{\font\=cmmi6\C}\kern.32mm\lower.1mm\hbox{$_{^{\emptyset,#1}}$}}% \|y\|\sNoreset_i produces \|y\|_{C^{\emptyset,i}}
\def\sbi#1{{_{\kern-0.1mm}}_{#1}} % same as _#1 but a little lower
\def\ais#1_#2{{}_{\font\=cmmi6\lower.15mm\hbox{\kern-.1mm\#1\kern.15mm}}\lower.3mm\hbox{${_{\kern-0.3mm_{#2}}}$}} %
\def\aais#1_#2{\kern.1mm{}_{\font\=cmmi6\lower.25mm\hbox{\kern-.1mm\#1\kern.15mm}}\lower.4mm\hbox{${_{\kern-0.3mm_{#2}}}$}} %
\def\ai#1{{}_{\font\=cmmi6\lower.15mm\hbox{\kern-.1mm\#1\kern.15mm}}} % 6:n pisteen kirjainalaindeksi
\def\yi#1{^{\font\=cmmi6\raise.0mm\hbox{\kern-.1mm\#1\kern.15mm}}} % 6:n pisteen kirjainyl"indeksi
\def\ar#1{{}_{\font\=cmr6\lower.15mm\hbox{\kern.1mm\#1}}} % 6:n pisteen numeroalaindeksi
\def\aar#1{_{\font\=cmr6\lower.15mm\hbox{\kern.1mm\#1}}} % 6:n pisteen numeroalaindeksi
\def\yr#1{^{\font\=cmr6\raise.0mm\hbox{\kern.3mm\#1}}} % 6:n pisteen numeroyl"indeksi
\def\yrai^#1_#2{^{\kern.4mm\hbox{\font\=cmr6\{#1}}}_{\kern.2mm{#2}}}
\def\upparentes#1{^{\kern.2mm\raise.2mm\hbox{\font\=cmr6\\char'050}\kern.1mm{#1}\kern.1mm\raise.2mm\hbox{\font\=cmr6\\char'051}}} % e.g. x\upparentes l gives x^{(l)} with 6point parentheses
\def\lupar{\kern.2mm\lower1mm\hbox{$^{^(}$}} % pieni vasen yl"sulku
\def\rupar{\lower1mm\hbox{$^{^)}$}\kern-.15mm} % pieni oikea yl"sulku
\def\yyi#1{^{\font\=cmmi6\lower.6mm\hbox{\kern-.25mm\#1\kern-.05mm}}} % pieni kirjainyl"indeksi, tarkoitettu derivaatan merkint""n
\def\yyr#1{^{\font\=cmr6\lower.45mm\hbox{\kern-.25mm\#1\kern-.15mm}}} % pieni numeroyl"indeksi, tarkoitettu derivaatan merkint""n
\def\yplus{\lower1mm\hbox{$^{^+}$}} % + merkki edelliseen
\def\yminus{\lower1mm\hbox{$^{^-}$}} % - merkki edelliseen
\def\aminus{{\kern.15mm\raise.3mm\hbox{$_{_-}$}\kern-.1mm}}%
\def\adot{\kern.2mm\hbox{\font\=cmb10\\char'056}}%
\def\ydot{\kern.2mm\raise1.9mm\hbox{\font\=cmb10\\char'056}}% k\ydot = k^. = the integer corresponding to the natural number k
\def\yydot{\kern.2mm\raise1.35mm\hbox{\font\=cmb7\\char'056}\kern.2mm}% the above for 7 point mode
\def\yydott{\kern.2mm\raise1.35mm\hbox{\font\=cmb6\\char'056}\kern.2mm}% the above for 7 point mode
\def\ClT{{\rm Cl}\kern.25mm\lower.4mm\hbox{$_{\Cal T}$}\kern0.2mm} % Cl_Cal T A = closure of A w.r.t. T
\def\IntT{\sp{\rm Int}\kern.2mm\lower.4mm\hbox{$_{\Cal T}$}\kern0.2mm} % Int_Cal T A = interior of A w.r.t. T
\def\Cl_taurd#1{\roman{Cl_{}}_{\kern0.37mm\hbox{\font\=cmmi8\\char'034}\kern-0.15mm{_{}}_{rd}\kern0.2mm#1\,}}% Cl_{tau_{rd}#1}S = closure of S in the topology of the topologized vector space E (= #1)

\def\iinc{\supseteq}
\def\exi#1{\exists\,#1\kern.2mm\,;}
\def\all#1{\forall\,#1\kern.2mm\,;}
\def\imply{\Rightarrow}

\def\spp{\kern0.07mm} % a horizontal very small positive space
\def\sp{\kern0.15mm} % a horizontal small positive space
\def\ssp{\kern0.37mm} % a bigger positive space
\def\snn{\kern-0.2mm} % a very small horizontal negative space
\def\sn{\kern-0.3mm} % a horizontal small negative space
\def\ssn{\kern-0.63mm} % a bigger negative space
\def\biggerlineskip#1 {\linebreak\nopagebreak\vskip-4.2mm\vskip.#1mm\nopagebreak\noindent}%
\def\Biggerlineskip#1 {\linebreak\nopagebreak\vskip-4.2mm\vskip#1pt\nopagebreak\noindent}%
\def\KP#1{\kern#1mm} % posive kern of #1 mm
\def\KN#1{\kern-#1mm} % negave kern of #1 mm
\def\nhskip#1mm{$\null$\kern#1mm}

\def\mhyppy#1{\null\kern#1mm}

\def\text#1{\hbox{\rm#1}}

\def\VBOX/#1/#2/HEREend{\vbox{#2\vskip-#1mm}\vfill\null\eject}
\def\œ$#1${\hbox{$#1$}} % text math which is not compressed or stretched
\def\"{\"a} \def\"{\"o}
\def\q#1{``\kern0.37mm#1\kern0.37mm"}
\def\newProCla#1\par#2\par{\vskip1.7mm\noindent\bf#1\it#2\vskip1.7mm}
\def\Prooff{{\font\=cmssi10\P\kern.37mmr\kern.37mmo\kern.37mmo\kern.37mmf\kern.37mm. }\rm}

\def\QED{\hfill\hbox{$\ \sqcap\kern-2.45mm\sqcup$}}

\def\noin{\noindent}
\def\Newline{\kern-10mm\newline}

\def\eps{\varepsilon}

\begin{document}

\title[$    \text{\sc On Nash\,--\,Moser applications}$]%
          {On an assertion about Nash\,--\,Moser applications}

\author[S. Hiltunen]{Seppo\ I\. Hiltunen}
\address{Helsinki University of Technology                             \vskip0mm$\hspace{2mm}$
           Institute of Mathematics, U311                              \vskip0mm$\hspace{2mm}$
           P.O.\ Box 1100                                              \vskip0mm$\hspace{2mm}$
           FIN-02015 HUT\vskip0mm
         FINLAND}
\email{shiltune\,@\,cc.hut.fi}



\begin{abstract}

By an example we show that Olaf M\"uller's assertion about his new theorems
being able to give anew some classical results previously obtained via
applications of Nash\,--\,Moser type theorems is unfounded. We also give
another example indicating some limitations in possible applications of
related new inverse function theorems.                        \end{abstract}

\maketitle

% ----------------------------------------------------------------------------

\noin Below, a {\it space\ssp} will mean a complete metrizable real locally
convex space $E\,$, hence a Fr\'echet space. By a {\it P-norm\ssp} for $E$ we
mean any function $\varrho:E\owns x\mapsto\varrho(x)=\|x\|$ such that $(x,y)
\mapsto\|x-y\|$ is a metric defining the topology of $E\,$.

\begin{definition}

A directionally differentiable (see \cite[Section 3]{Hic}\sp) map $\tilde f=
(E,F,f)$ of spaces $E,F$ where $f$ is a function $E\iinc\dom f\to F$, we call
{\it almost M-tame\ssp} at $x$ iff there are P-norms $\varrho_1$ for $E$ and $
\varrho_2$ for $F$ such that for any $z,u$ with $\varrho_1(z)\le 1$ and $u\in
E$ we have $x+z\in\dom f$ and $\varrho_2(\delta f(x+z,u)-\delta f(x,u))\le
\varrho_1(u)\,$.

\end{definition}

\begin{example}\label{Mul-cou}

Let $E$ be the Fr\'echet space of 1--periodic smooth functions $x:\Re\to\Re\,
$, and let the smooth $\varphi:\Re\to\Re$ be 1--periodic. Letting $\Iota:t
\mapsto n\,t$ and $\bold n:t\mapsto n$ for some fixed $n\in\Ze\setminus\{0\}\,
$, we consider the map $\tilde f\subtext{Ex\,2}=(E,E,f)$ where $f$
is defined by $U\owns x\mapsto\varphi\circ(\Iota+x)\cdot(\bold n+x')$ with $U$
being the set of $x$ in $E$ having $
    0 < \inf\ssp\{\,|\,n+x'(s)\ssp|:s\in\Re\sp\,\}\,$.         \end{example}

\begin{proposition}\label{xx}

If $\tilde f\subtext{Ex\,2}$ is almost M-tame at any $x\,$, then $\rng\varphi$
is a singleton.                                            \end{proposition}

\begin{proof} To proceed by reductio ad absurdum, we assume that $
\tilde f\subtext{Ex\,2}$ is almost M-tame at some $x$ and that $\varphi$ is
not constant, and we derive a contradiction. Now, there is some $t_0$ such
that $\varphi'(t_0)\not=0\,$, and we can also choose $s_0$ so that $
n\,s_0+x(s_0)=t_0\,$. There also are P-norms $u\mapsto\|u\|_1$ and $
v\mapsto\|v\|_2$ for $E$ such that for $\|z\|_1\le 1$ and for $u\in E$ we have
$\|\delta f(x+z,u)-\delta f(x,u)\|_2\le\|u\|_1\,$.

The topology of $E$ also being defined by the sequence of norms $\seq{p_i:i\in
\No}$ where $p_i(u)=\sup\{|u^{(l)}(s)|:s\in\Re$ and $l\le i\}\,$, there is an
odd $k\in\N$ such that $p_{k-1}(z)\le k^{-1}\imply\|z\|_1\le 1\,$. Then there
is $\delta_1\in\Rep$ such that $\|v\|_2\le\delta_1\imply p_{k-1}(v)\le 1\,$,
and there is $l\in\N$ such that $p_l(u)\le l^{-1}\imply\|u\|_1\le\delta_1\,$.

Taking $\eps\ar 0=l^{-1}$ and $m\in\N$ to be fixed below, we let $u:s\mapsto
\eps_0 $ and $z:s\mapsto(2\pi\,m)^{-k+\frac12}\sin(m\,2\pi\,(s-s_0))\,$. Then
putting $v=\delta f(x+z,u)-\delta f(x,u)\,$, we have $\|v\|_2\le\|u\|_1$ if $
p_{k-1}(z)\le k^{-1}\,$, and

$\mhyppy{11} v = \varphi'\circ(\Iota+x+z)\cdot u\cdot(\bold n+x'+z')
                + \varphi\circ(\Iota+x+z)\cdot u'$ 

$\mhyppy{33}    - \varphi'\circ(\Iota+x)\cdot u\cdot(\bold n+x')
                - \varphi\circ(\Iota+x)\cdot u'$

\centerline{$\phantom v = \eps_0(\varphi'\circ(\Iota+x+z)\cdot(\bold n+x'+z')
                         - \varphi'\circ(\Iota+x)\cdot(\bold n+x'))\,$,}

\noin and hence $v^{(k-1)}=\eps_0(\varphi'\circ(\Iota+x+z)\cdot z^{(k)}+\roman T_z)$
where $\roman T_z$ denotes a sum of terms where only derivatives
$z^{(i)}$ of $z$ with $i<k$ occur. Consequently, there is $M\in\Rep$ with
[ $p_{k-1}(z)\le k^{-1}$ and $s\in\Re\imply |\roman T_z(s)|\le M$ ] .

Now fixing $m$ so that $(2\pi\,m)^{-\frac12}\le k^{-1}$ and $
l+M < (2\pi\,m)^{\frac12}|\varphi'(t_0)|\,$, we have $p_{k-1}(z)\le k^{-1}$
and $p_l(u)\le l^{-1}\,$, hence $\|v\|_2\le\|u\|_1\le\delta_1\,$, and
consequently $p_{k-1}(v)\le 1<\eps_0((2\pi\,m)^{\frac12}|\varphi'(t_0)|-M)\le
\eps_0|\varphi'(n\,s_0+x(s_0))\cdot z^{(k)}(s_0)+\roman T_z(s_0)|
= |v^{(k-1)}(s_0)|\le p_{k-1}(v)\,$, a contradiction.            \end{proof}

In \cite[Section 5, Applications 1\,--\,5, pp.\ 25\,--\,26]{Mul} it is claimed
that certain previous results in \cite{Ha} can be obtained anew by the
approach in \cite{Mul}\ssp. In four of the asserted applications it is needed
to know that for a fixed tensor field $\omega$ on a manifold $N$ the pullback
map $\chi\mapsto\chi^*\omega$ by immersions $\chi:M\to N$ is
\q{bounded-differentiable}. It is obvious that if this indeed is the case,
then any local representation of the pullback map also is everywhere almost
M-tame.

It follows from our Proposition \ref{xx} above that this generally fails in
the case where $M=N=\mathbb S^1$ since our space $E$ is both a model space for
the manifold of immersions and linearly homeomorphic to the space of smooth
1--forms. Namely, in a local representation the pullback map for 1--forms has
precisely the form $x\mapsto\varphi\circ(\Iota+x)\cdot(\bold n+x')\,$. Similar
results can be obtained for other pairs of manifolds $M\sp,N$ and types of
tensors relevant to \cite[pp.\ 25\,--\,26]{Mul} by suitably adapting the idea
in the proof of Proposition \ref{xx} above.

\begin{example}\label{Hel-cou}

Let $E$ be the Fr\'echet space of smooth functions $x:[0,1]=I\to\Re\,$, and
let $\varphi:\Re\to\Re$ be a smooth diffeomorphism. We consider the
diffeomorphism $\tilde f\subtext{Ex\,4}=(E,E,f)$ where $f$ is defined by $
x\mapsto\varphi\circ x\,$.                                     \end{example}

\begin{proposition}\label{Helge B}

If $\tilde f\subtext{Ex\,4}$ is almost M-tame at $x\,$, then $\varphi''(t)=0$
for all $t\in\rng x\,$.                                    \end{proposition}

\begin{proof} Letting $\tilde f\subtext{Ex\,4}$ be almost M-tame at $x\,$, we
assume that $\varphi''(t_0)\not=0$ for some $t_0\in\rng x\,$, and show that a
contradiction follows. Now, there is some $s_0\in I$ with $x(s_0)=t_0\,$.
There also are P-norms $u\mapsto\|u\|_1$ and $v\mapsto\|v\|_2$ for $E$ such
that for $\|z\|_1\le 1$ and for $u\in E$ we have $
           \|\delta f(x+z,u)-\delta f(x,u)\|_2\le\|u\|_1\,$.

The topology of $E$ also being defined by the sequence of norms \œ$\seq{\,\sp
p\sbi i:i\in\No\ssp}$ where $p_i(u)=\sup\{|u^{(l)}(s)|:s\in\Re$ and $l\le i\}\,
$, there is an
odd $k\in\N$ such that $p_{k-1}(z)\le k^{-1}\imply\|z\|_1\le 1\,$. Then there
is $\delta_1\in\Rep$ such that $\|v\|_2\le\delta_1\imply p_k(v)\le 1\,$, and
there is $l\in\N$ such that $p_l(u)\le l^{-1}\imply\|u\|_1\le\delta_1\,$.

Taking $\eps\ar 0=l^{-1}$ and $m\in\N$ to be fixed below, we let $u=I\times
\{\eps_0\}$ and \œ$z=\Seq{\,(\sp 2\,\pi\,m\sp)^{\sp\mminus k\ssp+
\ssp\frac12}\ssp\sin\,(\sp 2\,\pi\,m\,(\sp s-s\ar 0)):s\in I\sp\,}\,$. Putting
$v=\delta f(x+z,u)-\delta f(x,u)\,$, the implication $p_{k-1}(z)\le k^{-1}
\imply\|v\|_2\le\|u\|_1$ holds, and we have\vskip.4mm

\centerline{$v = \varphi'\circ(x+z)\cdot u - \varphi'\circ x\cdot u
               = \eps_0(\varphi'\circ(x+z) - \varphi'\circ x)\,$,}\vskip.4mm

\noin and hence $v^{(k)}=\eps_0(\varphi''\circ(x+z)\cdot z^{(k)}+\roman T_z)$
where $\roman T_z$ denotes a sum of terms where only derivatives
$z^{(i)}$ of $z$ with $i<k$ occur. Consequently, there is $M\in\Rep$ such that
the implication [ $p_{k-1}(z)\le k^{-1}$ and $s\in I\imply |\roman T_z(s)|\le M$ ]
holds.

Now fixing $m>(2\pi)^{-1}\cdot\max\,\{k^2,|\varphi''(t_0)|^{-2}(l+M)^2\}\,$,
we have $p_{k-1}(z)\le k^{-1}$ and $p_l(u)\le l^{-1}\,$,
hence $\|v\|_2\le\|u\|_1\le\delta_1\,$, and consequently
$p_k(v)\le 1<\eps_0((2\pi\,m)^{\frac12}|\varphi''(t_0)|-M)\le
\eps_0|\varphi''((x+z)(s_0))\cdot z^{(k)}(s_0)
+\roman T_z(s_0)|=|v^{(k)}(s_0)|\le p_k(v)\,$, a contradiction. \end{proof}

From Proposition \ref{Helge B} it follows that if $\tilde f\subtext{Ex\,4}$
satisfies the premise in \cite[Theorem B, p.\ 3]{Gl} at $x\,$, then $
\varphi\,|\,\rng x$ is affine. Furthermore, if there is a set $A$ of functions
$x$ such that at every $x\in A$ the premise holds, and if in addition the set
$\,\rng\sn\bigcup\ssp A=\{\,x\ssp(t):x\in A$ and $t\in I\sp\,\}\,$ is dense in
the real line, then $\varphi$ itself is affine.

Noting that the idea in the proof of Proposition \ref{Helge B} can be
naturally adapted to show more generally that for a fixed smooth \œ$\varphi:
\varPi\spp\sqcap\varPi\aar 1\iinc\Omega\,|\,Q\to\varPi\aar 2$ only locally
affine maps \œ$\sp\tilde f\subtext{Re\,6}=(\sp\Cinfty(\spp Q\sbi{\varPi}\ssp,
\varPi\aar 1)\ssp,\Cinfty(\spp Q\sbi{\varPi}\ssp,\varPi\aar 2)\ssp,f\sp)$ with
$f$ given by \œ$x\mapsto$ \œ$\varphi\circ[\,\sp\roman{id\,},x\,]$ can
generically satisfy the premise in \cite[Theorem B, p.\ 3]{Gl}\ssp, we thus
see that problems involving this kind of maps $\sp\tilde f\subtext{Re\,6}$ are
practically ruled outside the domain of possible applications.

\begin{remark}

Our preceding observation, of course, does not preclude the possibility that
the theorems in the carefully written exposition \cite{Gl} may have nontrivial
applications where maps like $\sp\tilde f\subtext{Re\,6}$ do not play an
essential role.

There also are some historical implicit and inverse function theorems beyond
Banach spaces sharing this property of having premises too strong in order for
these theorems to be able to have serious applications involving maps of the
type $\sp\tilde f\subtext{Re\,6}$ above. We shall soon publish one such result
concerning Seip's various inverse and implicit function theorems in
\cite{Se72}\ssp.                                                \end{remark}

% ----------------------------------------------------------------------------

\end{document}